\documentclass[12pt]{article}
\usepackage[a4paper]{geometry}
\usepackage{amsfonts,amssymb,amsmath,enumerate}
\usepackage{color}
\usepackage{theorem,cite}
\usepackage{indentfirst}
\usepackage{textcomp}
\usepackage{mathdots}
\usepackage{bm}

\newtheorem{thm}{Theorem}[section]

\newtheorem{lem}[thm]{Lemma}
\newtheorem{coro}[thm]{Corollary}
\newtheorem{defi}[thm]{Definition}
\newcommand{\prf}{\textsl{Proof:}\ }
\newcommand{\qed}{\null \hfill {\rule{5pt}{5pt}}\\ \indent}

\theorembodyfont{\rmfamily}
\newtheorem{rmk}{Remark}[section]
\newcommand{\mb}[1]{\hspace{2.1ex}\mbox{#1}\hspace{2.1ex}}
\newcommand{\nonu}{\nonumber\\ }

\def\half{\frac12}

\newcommand\wh[1]{\widehat{#1}}

        \def\cU{{\cal U}}

\def\fe{{\mathfrak e}}
\def\ff{{\mathfrak f}}

\def\fh{{\mathfrak h}}

\newcommand{\CC}{{\mathbb C}}

\newcommand{\II}{{\mathbb I}}

\newcommand{\ZZ}{{\mathbb Z}}

          \def\sL{\mathsf{L}}
          
          \def\sR{\mathsf{R}}

\def\sv{\mathsf{v}}
\def\sw{\mathsf{w}}

\newcommand\Aqp[1]{\mathcal{A}_{q,p}(\mathfrak{\widehat{gl}}_{#1})}
\newcommand\Uq[1]{\mathcal{U}_q(\mathfrak{\widehat{gl}}_{#1})}

\DeclareMathOperator{\im}{im}
\DeclareMathOperator{\id}{Id}

\numberwithin{equation}{section}

\begin{document}
\pagestyle{empty}
\setcounter{page}{0}
\rightline{LAPTh-008/18}
\rightline{March 2018}
\vfill

\begin{center}

 {\LARGE  {\sffamily The quantum determinant \\ of the elliptic algebra $\mathcal{A}_{q,p}(\mathfrak{\widehat{gl}}_N)$} }\\[1cm]

\vspace{10mm}
  
{\Large 
 L. Frappat\footnote{luc.frappat@lapth.cnrs.fr},
   D. Issing\footnote{daniel.issing@lapth.cnrs.fr}
and E. Ragoucy\footnote{eric.ragoucy@lapth.cnrs.fr}
}\\[.41cm] 
Laboratoire d'Annecy-le-Vieux de Physique Th{\'e}orique LAPTh\footnote{Address: BP 110 Annecy-le-Vieux, F-74941 Annecy Cedex, France.} \\
Univ. Grenoble Alpes, USMB, CNRS, F-74000 Annecy
\end{center}
\vfill

\begin{abstract}
We introduce the quantum determinant for the elliptic quantum algebra $\mathcal{A}_{q,p}(\mathfrak{\widehat{gl}}_N)$ and prove that it generates the center of this algebra. We also show that it is group-like for the quasi-Hopf structure, which allows us to define the 
elliptic quantum algebra $\mathcal{A}_{q,p}(\mathfrak{\widehat{sl}}_N)$.
\end{abstract}

\vfill
\newpage
\pagestyle{plain}
\section{Introduction} 

It is well known that certain algebraic structures are characterized by solutions of the Yang--Baxter equation (YBE) with spectral parameter: for example, the rational and trigonometric solutions lead to the Yangian and quantum affine algebras respectively. 
In this framework, called the FRT formalism \cite{FRT}, the generators of the algebra under consideration are encapsulated into the Lax matrix $L(z)$, whose intertwining properties are coded in the $R$-matrix (RLL relations). 
Within this approach, Foda \emph{et al.} \cite{FIJKMY,FIJKMY95} were able to define the elliptic quantum algebra $\Aqp{2}$ using the elliptic solution of the YBE ($R$-matrix of the 8-vertex model) discovered by Baxter \cite{Baxter}. 
This solution exhibits entries expressed in terms of the Jacobi theta function and depends on a deformation parameter $q$ and an elliptic nome $p$. 

When properly normalized, Baxter's elliptic $R$-matrix satisfies $R_{12}(q^{-1})=A_{12}$ where $A_{12}$ is the antisymmetrizer on $(\CC^2)^{\otimes 2}$. 
This property implies a quantum determinant formula, such a quantum determinant lying in the center of the algebra. 
This allows to define the elliptic algebra $\mathcal{A}_{q,p}(\mathfrak{\widehat{sl}}_{2})$ by imposing further a restriction on the quantum determinant, namely
\begin{equation}
\mathcal{A}_{q,p}(\mathfrak{\widehat{sl}}_{2}) = \mathcal{A}_{q,p}(\mathfrak{\widehat{gl}}_{2}) / 
\langle \text{qdet}\,L(z) - q^{c/2} \rangle
\label{eq:Aqpsl2}
\end{equation}
where $c$ is the central charge of the algebra.

The generalization to the $\mathfrak{\widehat{gl}}_{N}$ case builds on the elliptic $\ZZ_N$-symmetric $R$-matrix derived by Belavin \cite{Bela}. 
This matrix reduces to Baxter's elliptic $R$-matrix when $N=2$.
Along the same lines as above, an elliptic quantum algebra $\Aqp{N}$ was defined in \cite{AFRS99} by using the elliptic $\ZZ_N$-symmetric $R$-matrix.
However, an explicit formula for the quantum determinant is still missing in the case of $\Aqp{N}$ when $N>2$. 
The main result of this paper is the derivation of the quantum determinant formula for $\Aqp{N}$ and the proof that it generates the center of the algebra.
Hence formula \eqref{eq:Aqpsl2} can be extended to generic $N$.

Let us remark that this derivation has a slightly different final step, compared to the usual approach used in the quantum affine or Yangian cases \cite{FRT,KS2,MNO}. 
Indeed, the usual approach uses  the fact that the antisymmetrizer on $(\CC^N)^{\otimes N}$ can be expressed as some
  product of $R$-matrices. 
In the elliptic case, this does not seem to be true when $N>2$. It does not prevent to apply the antisymmetrizer to an $N$-fold product of Lax matrices and
 to construct an analog of the quantum determinant that is still central. However, it remains to compute the value of the quantum determinant in the fundamental evaluation representation. Without the expression of the antisymmetrizer in term of $R$-matrices, it requires some more work with respect to the usual case.

The plan of the paper is as follows. 
In section \ref{sec:Aqp}, we recall the definition and basic properties of the elliptic quantum algebra $\Aqp{N}$. 
Section \ref{sec:UqN} is a reminder devoted to the quantum affine algebra $\Uq{N}$. 
Here, we note in particular that the limit $p \to 0$ of the elliptic $R$-matrix is a \emph{twisted} version of the affine $R$-matrix in the principal gradation. 
The main result of the paper, the formula for the quantum determinant of $\Aqp{N}$, is stated in section \ref{sec:main}. 
Proofs are gathered in section \ref{sec:proofs}.


\section{The elliptic quantum algebra $\Aqp{N}$\label{sec:Aqp}}

\subsection{Definition of $\Aqp{N}$}

Let us briefly review the construction of the elliptic quantum algebra $\Aqp{N}$.
We consider a free, associative, unital algebra generated by operators $L_{i,j}[n]$ ($1 \le i,j \le N$, $n \in \ZZ$), compactly represented by the formal series of the spectral parameter $z \in \CC$
\begin{equation}
L_{i,j}(z)=\sum_{n \in \mathbb{Z}} L_{i,j}[n] \, z^{-n}  
\label{seriesgenerators}
\end{equation}
and encapsulated into the so-called Lax matrix
\begin{equation*}
L(z)=\sum_{i,j=1}^N L_{i,j}(z) \, e_{i,j} \,,
\label{Lax}
\end{equation*}
where the matrices $e_{i,j}$ are the usual elementary matrices. 
We adjoin an invertible central element written as $q^c$ ($q$ is the deformation parameter and $c$ the central charge).

The $\Aqp{N}$ algebra is defined by imposing a set of exchange relations (commonly refereed to as RLL or FRT relations \cite{FRT}) on these generators:
\begin{equation}
\wh{R}_{12}(\frac{z}{w})L_1(z)L_2(w)=L_2(w)L_1(z)\wh{R}^{*}_{12}(\frac{z}{w}) \,.
\label{RLL}
\end{equation}
The indices refer to the spaces in which the Lax operators and the matrix $\wh{R}(z)$ operate, with $L_1(z)=L(z) \otimes \mathbb{I}$ and $L_2(z)=\mathbb{I}\otimes L(z)$.
In addition to its spectral parameter dependence, the $R$-matrix $\wh{R}(z)$  depends on two complex parameters $q$ (deformation parameter) and $p$ (elliptic nome). To lighten presentation, unless needed, we will not write explicitly the dependence in $q$ and $p$. The second $R$-matrix appearing in  \eqref{RLL} is related to the original one through $\wh{R}^{*}_{12}(z)=\wh{R}_{12}(z)\big|_{p\to p^{*}=p q^{-2c}}$. 

The $\Aqp{N}$ algebra is a quasi-Hopf algebra with coproduct
\begin{eqnarray}\label{coproduct}
\Delta L(z;p) &=& \Big(1{\otimes}L(zq^{c^{(1)}/2};p)\Big)\cdot \Big(L(zq^{-c^{(2)}/2};pq^{-2c^{(2)}}){\otimes}1\Big),\\
\Delta c &=& c\otimes 1 + 1\otimes c\equiv c^{(1)}+c^{(2)}
\end{eqnarray}
where we have specified the $p$ dependence in $L(z)$, due to shifts occurring when we apply $\Delta$ (for more details see \cite{JKOS,Konno}). One can easily check that $\Delta$ is an algebra morphism.

In order to define explicitly the matrix $\wh{R}(z)$ appearing in here, we  first introduce a slightly modified matrix $R(z)=\sum R_{a,c}^{b,d}(z) \, e_{a,b} \otimes e_{c,d}$, whose non-vanishing entries obey $a+c=b+d$, the addition of indices being understood modulo $N$. 
For $1 \le a,b,c,d \le N$, one defines
\begin{equation}
R_{a,c}^{b,d}(z) = \eta (z) S_{a,c}^b(z) \, \omega^{(a+c-b-d)/2} \, \delta_{a+c,b+d}^{(\text{mod}\;N)} \;,
\end{equation}
where
\begin{equation}
S_{a,c}^b(z) = z^{\frac{2(b-a)}{N}} q^{\frac{2(c-b)}{N}} p^{\frac{(b-a)(c-b)}{N}}\frac{\Theta_{p^N}(p^{N+c-a} q^2 z^2)}{\Theta_{p^N}(p^{N+c-b} z^2)\Theta_{p^N}(p^{N+b-a} q^2)} \;.
\label{RMatrixAN}
\end{equation}
Remark that since $\omega=e^{2i\pi/N}$, the factor $\omega^{(a+c-b-d)/2}$ equals $\pm 1$. This factor results from the gauge transformation done on the Belavin $R$-matrix that allows one to recover the original $\Aqp{2}$ matrix. Note that in \eqref{RMatrixAN} one can consider any index $c \in \ZZ$ because of the identity $S_{a,c+N}^b(z) = S_{a,c}^b(z)$. 

The normalization coefficient $\eta(z)$ is given by
\begin{equation}
\eta(z) = \frac{z^{\frac{2}{N}}}{\kappa_N(z^2)} \frac{(p^N,p^N)^3_{\infty}}{(p,p)^3_{\infty}}\frac{\Theta_{p}(q^2)\Theta_{p}(p z^2)}{\Theta_{p}(q^2 z^2)}
\end{equation}
with 
\begin{equation}
\frac{1}{\kappa_N(z^2)} = 
\frac{(q^{2N} z^{-2};p,q^{2N})_{\infty} (q^{2} z^{2};p,q^{2N})_{\infty} (p z^{-2};p,q^{2N})_{\infty} 
(p q^{2N-2} z^{2};p,q^{2N})_{\infty}}{(q^{2N} z^{2};p,q^{2N})_{\infty} (q^{2} z^{-2};p,q^{2N})_{\infty} 
(p z^{2};p,q^{2N})_{\infty} (p q^{2N-2} z^{-2};p,q^{2N})_{\infty}} \;.
\label{kappa}
\end{equation}
The matrix $S$ is $\mathbb{Z}_N$-symmetric by construction, i.e. the coefficients satisfy $S_{a+n,c+n}^{b+n}(z)=S_{a,c}^b(z)$, $n=1,...,N$, where again the addition of indices is modulo $N$. 

Here, $\Theta_p(z)$ denotes the Jacobi theta function defined by ($p \in \CC$ such that $|p|<1$)
\begin{equation}\label{theta:prod}
\Theta_p(z) = (z;p)_\infty \, (pz^{-1};p)_\infty \, (p;p)_\infty 
\end{equation}
and the infinite $q$-Pochhammer symbols are given by 
\begin{equation}
(z;p_1,\dots,p_m)_\infty = \prod_{n_i \ge 0} (1-z p_1^{n_1} \dots p_m^{n_m}) \,.
\end{equation}
It is easy to show that the Jacobi theta function enjoys the following properties:
\begin{align}\label{id:theta}
& \Theta_{a}(az) = \Theta_{a}(z^{-1}) 
\mb{and} 
\Theta_{a}(a^nz) = \frac{(-1)^n}{z^n\,a^{n(n-1)/2}}\,\Theta_{a}(z)\,,\quad \forall\ n\in\ZZ,  \\
& \Theta_{a^2}(az) = \Theta_{a^2}(az^{-1}).
\label{id:theta2}
\end{align}

The matrices $\wh{R}(z)$ and $R(z)$ differ only by a suitable normalization, which reads
\begin{equation}
\wh{R}(z) = \tau_N(q^{\frac{1}{2}} z^{-1}) R(z) \,,
\label{RealR}
\end{equation}
where
\begin{equation}
\tau_N(z)=z^{\frac{2}{N}-2} \frac{\Theta_{q^{2N}}(q z^2)}{\Theta_{q^{2N}}(q z^{-2})} \,.
\label{rescale}
\end{equation}

\subsubsection*{Example: the $\Aqp{2}$ algebra}

Specifying $N=2$, the $R$-matrix \eqref{RMatrixAN} reads explicitly as:
\begin{equation}
R(z)=\frac{1}{\kappa_2(z^2)}\frac{(p^2;p^2)_{\infty}}{(p;p)^2_{\infty}}
\begin{pmatrix}
a(z) & 0 & 0 & d(z)\\
0 & b(z) & c(z) & 0\\
0 & c(z) & b(z) & 0\\
d(z) & 0 & 0 & a(z)
\end{pmatrix}
\label{RMatrixA2}
\end{equation}
and the entries take the following form:
\begin{eqnarray}
\begin{aligned}[l]
a(z) &=z^{-1} \frac{\Theta_{p^2}(p z^2)\Theta_{p^2}(p q^2)}{\Theta_{p^2}(p q^2 z^2)} \,, \\
b(z)&=q z^{-1} \frac{\Theta_{p^2}( z^2)\Theta_{p^2}(p q^2)}{\Theta_{p^2}(q^2 z^2)} \,,
\end{aligned}
\hspace{1cm}
\begin{aligned}[l]
c(z) &= \frac{\Theta_{p^2}(p z^2)\Theta_{p^2}(q^2)}{\Theta_{p^2}(q^2 z^2)} \,, \\
d(z)&=-\frac{p^{\frac{1}{2}}}{q z^2}\frac{\Theta_{p^2}(z^2)\Theta_{p^2}(q^2)}{\Theta_{p^2}(p q^2 z^2)} \,.
\label{EntriesA3}
\end{aligned}
\end{eqnarray}
The $\ZZ_2$-symmetry of the $R$-matrix corresponds to a symmetry w.r.t. the diagonal  and w.r.t. the anti-diagonal. 
It implies in particular that $R(z)^{t_1t_2}=R(z)$. This is no longer the case for higher rank algebras.

\begin{rmk}
It has been shown \cite{Fron,JKOS} that the elliptic quantum algebras are quasi-Hopf algebras, obtained by a twisting procedure on the quantum affine algebra $\Uq{N}$. The explicit expression of the Drinfel'd twist and the proof that this twist indeed satisfies the shifted cocycle condition were given in \cite{JKOS}, see also \cite{ABRR}. In the case of $\Aqp{2}$, the evaluated $R$-matrix \eqref{RMatrixA2} was recovered through this procedure. Note that for $N>2$, the explicit expression of the twist in the $N$-dimensional evaluation representation is, to our knowledge, still unknown.
\end{rmk}

\subsection{Properties of the $R$-matrix}

Let $g$ and $h$ be the matrices of order $N$ defined by $g_{ij} = \omega^i\,\delta_{ij}$ and $h_{ij} = \delta_{i+1,j}$ for $1\leq i,j\leq N$ with $\omega = e^{2i\pi/N}$, the addition of indices being understood modulo $N$.

We recall the following properties of the $R$-matrix $R(z)$ \cite{Tra85,RT86}:
\begin{itemize}
\item Yang--Baxter equation (also holds for $\wh{R}(z)$):
\begin{equation}
R_{12}(\frac{z_1}{z_2}) \, R_{13}(\frac{z_1}{z_3}) \, R_{23}(\frac{z_2}{z_3})=
R_{23}(\frac{z_2}{z_3}) \, R_{13}(\frac{z_1}{z_3}) \, R_{12}(\frac{z_1}{z_2}) \,.
\label{YBE}
\end{equation}
\item Unitarity:
\begin{equation}
R_{12}(z) \, R_{21}(z^{-1}) = \II \,,
\label{eq222}
\end{equation}
\item Regularity ($P_{12}$ is the permutation matrix):
\begin{equation}
R_{12}(1) = P_{12} \,,
\end{equation}
\item Crossing-symmetry:
\begin{equation}
R_{12}(z)^{t_2} \, R_{21}(z^{-1}q^{-N})^{t_2} = \II \,,
\label{eq223}
\end{equation}
\item Antisymmetry:
\begin{equation}
R_{12}(-z) = \omega \, (g^{-1} \otimes \II) \, R_{12}(z) \, (g \otimes \II) \,,
\label{eq224}
\end{equation}
\item Quasi-periodicity:
\begin{equation}
\widehat R_{12}(-z p^{\frac{1}{2}}) = (g^{\frac{1}{2}} h g^{\frac{1}{2}} \otimes \II)^{-1} \, \widehat R_{21}(z^{-1})^{-1} \, (g^{\frac{1}{2}} h g^{\frac{1}{2}} \otimes \II) \,,
\label{eq225}
\end{equation}
\item Invariance:
\begin{equation}
(h \otimes h) \, R_{12}(z) =  R_{12}(z) \, (h \otimes h) \,.
\label{eq:h-inv}
\end{equation}
\end{itemize}

\begin{rmk}
The crossing-symmetry and the unitarity properties of $R_{12}$ allow to exchange the inversion and the
transposition when applied to the matrix $R_{12}$ (or to the matrix $\widehat R_{12}$). It provides a \textsl{crossing-unitarity relation} (also valid for $\widehat R$ thanks to the $q^N$-periodicity of the function $\tau_N$):
\begin{equation}
\Big(R_{12}(x)^{t_2}\Big)^{-1} = \Big(R_{12}(q^Nx)^{-1}\Big)^{t_2} \,.
\label{eq230}
\end{equation}
Note also that the unitarity property for $\widehat R_{12}$  reads 
\begin{equation}
\widehat R_{12}(z) \, \widehat R_{21}(z^{-1}) = \tau_N(q^{\frac 12}z) \, \tau_N(q^{\frac 12}z^{-1}) \equiv \mathcal{U}(z),
\label{eq:unitarity}
\end{equation}
where the function $\cU(z)$  is defined as
\begin{equation}\label{def:U}
\cU(z)=q^{\frac2N-2}\,\frac{\Theta_{q^{2N}}(q^2z^2) \, \Theta_{q^{2N}}(q^2z^{-2})} {\Theta_{q^{2N}}(z^2)\Theta_{q^{2N}}(z^{-2})} .
\end{equation}
\end{rmk}

\section{The quantum affine algebra $\Uq{N}$\label{sec:UqN}}

We recall in this section the different gradations that may be used in constructing the $R$-matrix defining the algebra $\Uq{N}$ in the FRT formalism. As will be shown below, the $R$-matrix obtained as the non-elliptic limit of \eqref{RealR} appears to be a twisted version of the $R$-matrix in the principal gradation.
 
\subsection{Homogeneous gradation.} 
In the FRT formalism, the quantum affine algebra $\Uq{N}$ is described as an associative algebra defined by generators and relations. 
The generators $L_{i,j}^\pm[\mp n]$, where $n\in\ZZ_{\ge 0}$, $1 \le i,j \le N$ and $L_{i,j}^+[0] = L_{j,i}^-[0] = 0$ for $i>j$, are coded in formal generating functions $L_{i,j}^\pm(z)$, themselves encapsulated into the Lax matrices $L^{\pm}(z)$:
\begin{equation}
L^{\pm}(z)=\sum_{i,j=1}^N L_{i,j}^{\pm}(z) \, e_{i,j}
\quad \text{and} \quad
L_{i,j}^{\pm}(z) = \sum_{n=0}^{\infty} L_{i,j}^{\pm}[\mp n] \, z^{\pm n} \,.
\label{eq:Laxhomog}
\end{equation}
The relations are the well-known RLL relations 
\begin{align}
\label{eq:RLLUq}
R_{12}(\frac{z_\pm}{w_\pm}) L_1^{\pm}(z) L_2^{\pm}(w) &= L_2^{\pm}(w) L_1^{\pm}(z) R_{12}(\frac{z_\pm}{w_\pm}) \,, \\
R_{12}(\frac{z_\pm}{w_\mp}) L_1^{+}(z) L_2^{-}(w) &= L_2^{-}(w) L_1^{+}(z) R_{12}(\frac{z_\mp}{w_\pm}) \,,
\end{align}
where $z_{\pm} = zq^{\pm c/2}$, $w_{\pm} = wq^{\pm c/2}$, $c$ is the central charge and the matrix $R_{12}(z)$ is given by
\begin{align}
R_{12}(z)=\rho_N(z)\left[\sum_i e_{i,i} \otimes e_{i,i} + \frac{q(1-z)}{1-q^2z} \sum_{i \ne j} e_{i,i} \otimes e_{j,j} + \frac{(1-q^2)}{1-q^2z} \left(\sum_{i<j}+z\sum_{i>j}\right)e_{i,j} \otimes e_{j,i}\right] \,.
\label{eq:Rhomog}
\end{align}
The normalization factor $\rho_N(z)$ expressed as
\begin{equation}
\rho_N(z)=q^{\frac{1}{N}-1} \frac{(q^2 z;q^{2N})_{\infty}(q^{2N-2} z;q^{2N})_{\infty}}{(z;q^{2N})_{\infty}(q^{2N} z;q^{2N})_{\infty}} \;.
\label{NormUq}
\end{equation}
Let $\fe_i,\ff_i$  ($0 \le i \le N-1$) and $\fh_i$ ($0 \le i \le N$) denote the generators of $\Uq{N}$ in the Serre--Chevalley basis and let $\mathcal{R}$ be the universal $R$-matrix of $\Uq{N}$ (see e.g. \cite{KT1994}).
The $R$-matrix \eqref{eq:Rhomog} is obtained from $\mathcal{R}$ by calculating its image $R(z/w) = (\pi_z \otimes \pi_w) \mathcal{R}$ in the $N$-dimensional evaluation representation $\pi_z$ such that ($1 \le i \le N$)
\begin{alignat}{3}
& \pi_z(\fe_i) = e_{i,i+1} \,, & \quad & \pi_z(\ff_i) = e_{i+1,i} \,, & \quad & \pi_z(\fh_i) = e_{i,i} \,, \\
& \pi_z(\fe_0) = z e_{N,1} \,, & \quad & \pi_z(\ff_0) = z^{-1} e_{1,N} \,, & \quad & \pi_z(\fh_0) = e_{N,N} - e_{1,1} \,.
\end{alignat}
This defines the so-called homogeneous gradation. \\
The quantum affine algebra $\Uq{N}$ is endowed with the following coproduct structure:
\begin{equation}
\Delta \big( L_{i,j}^{\pm}(z) \big) = \sum_{k=1}^N L_{k,j}^{\pm}(zq^{\mp c^{(2)}/2}) \otimes L_{i,k}^{\pm}(zq^{\pm c^{(1)}/2}) \;,
\end{equation}
where $c^{(1)} = c \otimes 1$ and $c^{(2)} = 1 \otimes c$. \\
The quantum determinant is given in the homogeneous gradation by
\begin{equation}\label{qdet-homo}
\text{qdet}\,L(z) = \sum_{\sigma \in \mathfrak{S}_N} \text{sgn}(\sigma) \, q^{\ell(\sigma)} L^{+}_{1,\sigma(1)}(z) \, L^{+}_{2,\sigma(2)}(zq^{-2}) \dots L^{+}_{N,\sigma(N)}(zq^{2-2N}) \;,
\end{equation}
where $\ell(\sigma)$ denotes the length of the permutation $\sigma$ and $\text{sgn}(\sigma)=(-1)^{\ell(\sigma)}$. \\
Finally, thanks to the RLL relations, the action of the finite Cartan generators on the Lax matrices is given by ($1 \le i,j,k \le N$)
\begin{equation}
q^{\fh_i} \, L^{\pm}_{j,k}(w) = L^{\pm}_{j,k}(w) \, q^{\fh_i+\delta_{ik}-\delta_{ij}} \;.
\label{eq:hL}
\end{equation}

\subsection{Principal gradation} 
Another possible choice is the principal gradation. In that case, the evaluation map $\pi_z$ is given by
\begin{alignat}{3}
& \pi_z(\fe_i) = z^{2/N} e_{i,i+1} \,, & \quad & \pi_z(\ff_i) = z^{-2/N} e_{i+1,i} \,, & \quad & \pi_z(\fh_i) = e_{i,i} \,, \\
& \pi_z(\fe_0) = z^{2/N} e_{N,1} \,, & \quad & \pi_z(\ff_0) = z^{-2/N} e_{1,N} \,, & \quad & \pi_z(\fh_0) = e_{N,N} - e_{1,1} \,.
\end{alignat}
The $R$-matrix in the principal gradation reads:
\begin{align}
\sR(z) &= \rho_N(z^2)\left[\sum_i e_{i,i} \otimes e_{i,i} + \frac{q(1-z^2)}{1-q^2z^2} e_{i,i} \otimes e_{j,j} \right. \nonumber \\
&\left. + \frac{z(1-q^2)}{1-q^2z^2} \left(\sum_{i<j}z^{(2j-2i-N)/N}+\sum_{i>j}z^{(2j-2i+N)/N}\right)e_{i,j} \otimes e_{j,i}\right] \,.
\label{eq:Rprinc}
\end{align}
The two matrices \eqref{eq:Rhomog} and \eqref{eq:Rprinc} are related by a gauge transformation 
\begin{equation}
\sR(z/w) = V(z) \otimes V(w) \, R(z^2/w^2) \, (V(z) \otimes V(w))^{-1}
\end{equation}
with $V(z) = \displaystyle\sum_{i=1}^{N} z^{(N+1-2i)/N}\,e_{i,i}$. \\
It follows that the Lax matrices $L^{\pm}(z)$ and $\sL^{\pm}(z)$ that define the quantum affine algebra in the homogeneous and principal gradations respectively are related by
\begin{align}
\label{eq:gaugeLp}
& \sL^+(z) = V(zq^{c/2}) \, L^+(z^2) \, V(zq^{-c/2})^{-1}  \,, \\
\label{eq:gaugeLm}
& \sL^-(z) = V(z) \, L^-(z^2) \, V(z)^{-1} \,.
\end{align}
Note that these relations ensure that equation \eqref{eq:hL} also holds for the Lax matrices $\sL^\pm(z)$ in the principal gradation. \\
The quantum determinant is then given in the principal gradation by
\begin{equation}\label{qdet-princ}
\text{qdet}\,\sL(z) = \sum_{\sigma \in \mathfrak{S}_N} \text{sgn}(\sigma) \, q^{\ell(\sigma) + \frac{2}{N} \sum_{i=1}^N i(\sigma(i)-i)} \, \sL^{+}_{1,\sigma(1)}(z) \, \sL^{+}_{2,\sigma(2)}(zq^{-1}) \dots \sL^{+}_{N,\sigma(N)}(zq^{1-N}) \,.
\end{equation}

\subsection{Non-elliptic presentation\label{sect:non-ell}} 
The limit $p \rightarrow 0$ of the $\Aqp{N}$ algebra allows us to reveal still another presentation of the quantum affine algebra $\Uq{N}$. Since this presentation is related to the non-elliptic limit of $\Aqp{N}$, we will call it the \textsl{non-elliptic presentation.}
The $R$-matrix obtained in this limit reads:
\begin{align}
\textsf{R}'(z) &= \rho_N(z^2)\left[\sum_i e_{i,i} \otimes e_{i,i} + \frac{q(1-z^2)}{1-q^2z^2} \left(\sum_{i<j}q^{(2j-2i-N)/N}+\sum_{i>j}q^{(2j-2i+N)/N}\right)e_{i,i} \otimes e_{j,j} \right. \nonumber \\
&\left. + \frac{z(1-q^2)}{1-q^2z^2} \left(\sum_{i<j}z^{(2j-2i-N)/N}+\sum_{i>j}z^{(2j-2i+N)/N}\right)e_{i,j} \otimes e_{j,i}\right] \,.
\label{eq:RUqell}
\end{align}
When $N>2$, this matrix differs from the previous one, essentially by some powers of $q$ in the diagonal terms. 
Obviously, it is still $\mathbb{Z}_N$-symmetric. 
It can be obtained from \eqref{eq:Rprinc} by a (constant non factorized) diagonal twist:
\begin{equation}
\textsf{R}'(z) = F_{21} \, \sR(z) \, F_{12}^{-1} 
\end{equation}
where
\begin{equation}
F_{12} = \sum_{i=1}^{N} e_{i,i} \otimes e_{i,i} + \sum_{1 \le i \ne j \le N} q^{\alpha_{ij}} \, e_{i,i} \otimes e_{j,j}
\end{equation}
with, for $i<j$, $\alpha_{ij}=\frac{1}{2}+(i-j)/N$ and $\alpha_{ji}=-\alpha_{ij}$. We set by convention 
$\alpha_{ii}=0$ for all $i$. Remark that for $N=2$, $\alpha_{12}=0$, so that the twist is $\II\otimes\II$.\\
The algebra is still defined by Eqs. \eqref{eq:RLLUq}--\eqref{eq:Laxhomog} where the Lax matrices $L{^\pm}(z)$ are now replaced by $\sL'{^\pm}(z)$. \\
At the universal level, the twisted $R$-matrix is given by 
\begin{equation}
\mathcal{R}^{\mathcal{F}} = \mathcal{F}_{21} \, \mathcal{R} \, \mathcal{F}_{12}^{-1}
\end{equation}
with
\begin{equation}
\mathcal{F}_{12} = q^{\sum_{ij} \alpha_{ij} \fh_{i} \otimes \fh_{j}} \,.
\label{eq:twistuniv}
\end{equation}
Here $\fh_i$ ($i=1,\dots,N$) are the Cartan generators of the finite quantum algebra $\mathcal{U}_q(\mathfrak{gl}_N)$ satisfying the following commutation relations ($j=1,\dots,N-1$):
\begin{equation}
[\fh_i,\fe_j] = (\delta_{ij}-\delta_{i,j+1}) \fe_j \;, \quad 
[\fh_i,\ff_j] = -(\delta_{ij}-\delta_{i,j+1}) \ff_j \;, \quad 
[\fe_j,\ff_j] = \dfrac{q^{\fh_j-\fh_{j+1}}-q^{\fh_{j+1}-\fh_j}}{q-q^{-1}} \,.
\end{equation}
The universal twist \eqref{eq:twistuniv} satisfies the cocycle condition $\mathcal{F}_{12} (\Delta \otimes \text{id}) \mathcal{F} = \mathcal{F}_{23} (\text{id} \otimes \Delta) \mathcal{F}$, ensuring that the universal $R$-matrix $\mathcal{R}^{\mathcal{F}}$ satisfies the Yang--Baxter equation while the $R$-matrix $\mathcal{R}$ does. \\
The relation between the corresponding Lax matrices $\sL^{\pm}$ and $\sL'^{\pm}$ can be expressed as
\begin{equation}
\sL'^{\pm}(z) = (\pi_z \otimes \text{id})\mathcal{F}_{21} \, \sL^{\pm}(z) \, (\pi_z \otimes \text{id})\mathcal{F}_{12}^{-1} \,.
\label{eq:corresL}
\end{equation}
In the evaluation representation $\pi_z$, one gets
\begin{align}
(\pi_z \otimes \text{id})\mathcal{F}_{12} = (\pi_z \otimes \text{id})\mathcal{F}_{21}^{-1} = \sum_{i=1}^{N} q^{\sum_{j=1}^N \alpha_{ij} \fh_j} \, e_{i,i} \,.
\end{align}
The twist being diagonal and depending only on the finite Cartan generators, the equation \eqref{eq:hL} also holds for the Lax matrices $\sL'^\pm(z)$. \\ 
The coproduct of the twisted algebra is given by $\Delta^{\mathcal{F}} = \mathcal{F}_{12} \, \Delta \, \mathcal{F}_{12}^{-1}$. A direct calculation shows that
\begin{equation}
\Delta^{\mathcal{F}} \big( \sL'^{\pm}_{i,j}(z) \big) = \sum_{k=1}^N  {\sL'^{\pm}_{k,j}}(zq^{\mp c^{(2)}/2}) \otimes {\sL'^{\pm}_{i,k}}(zq^{\pm c^{(1)}/2})
\end{equation}
from which it follows that the twisted algebra gets the same coproduct structure as the original algebra.
\\
Applying the twist to the expression \eqref{qdet-princ}, and using the correspondence \eqref{eq:corresL},
we get an expression for the quantum determinant in this new presentation. 
It is again expressed as a sum over permutations:
\begin{equation}
\text{qdet}\,\sL(z) = \sum_{\sigma \in \mathfrak{S}_N} \text{sgn}(\sigma) \, q^{n_\sigma} \, \sL'^{+}_{1,\sigma(1)}(z) \, \sL'^{+}_{2,\sigma(2)}(zq^{-1}) \dots \sL'^{+}_{N,\sigma(N)}(zq^{1-N}) \,,
\end{equation}
where $n_\sigma=\ell(\sigma) + \frac{2}{N} \sum_{i=1}^N i(\sigma(i)-i) + \sum_{1 \le i<j \le N} (\alpha_{\sigma(i),\sigma(j)} - \alpha_{ij}) $. A detailed analysis of $n_\sigma$, using the explicit expression of the coefficients $\alpha_{ij}$, shows that it vanishes identically. 
Then, the quantum determinant is given in the non-elliptic limit by
\begin{equation}\label{qdet-nonell}
\text{qdet}\,\sL'(z) = \sum_{\sigma \in \mathfrak{S}_N} \text{sgn}(\sigma) \, \sL'^{+}_{1,\sigma(1)}(z) \, \sL'^{+}_{2,\sigma(2)}(zq^{-1}) \dots \sL'^{+}_{N,\sigma(N)}(zq^{1-N})\,.
\end{equation}
Let us remark that the relation \eqref{qdet-nonell} is based on the (undeformed) antisymmetrizer, contrarily to the expressions found for the homogeneous and principal gradations that are based on $q$-deformed versions of it.
When $N=2$, $\sR(z)$ and $\sR'(z)$ coincide, and only the homogeneous gradation provides a deformed antisymmetrizer.

\section{Main results\label{sec:main}}

We are now in position to state the main result of this article.
\begin{thm}\label{thm:qdet}
Let $A^{(N)}_{N}$ be the antisymmetrizer of $N$ spaces $\CC^N$. For generic values of the parameters $p$, $q$ and of the central charge $c$, 
one has the following identity
\begin{equation}
L_1(z)\dots L_N(z q^{1-N})\, A^{(N)}_{N} = \mathrm{qdet} L(z)\, A^{(N)}_{N} \,,
\label{qdet0}
\end{equation}
where $\mathrm{qdet} L(z)$, called the  \textit{quantum determinant}, is a scalar function that lies in the center of the $\Aqp{N}$ algebra. It can be rewritten as
\begin{eqnarray}
\mathrm{qdet} L(z)&=&\mathrm{tr}_{1\dots N}\left(L_1(z)\dots L_N(z q^{1-N}) A^{(N)}_{N}\right)
\label{qdet2}
\\
&=& \sum_{\sigma \in \mathfrak{S}_N} \mathrm{sgn}(\sigma)L_{1,\sigma(1)}(z) L_{2,\sigma(2)}(\frac{z}{q})\dots L_{N,\sigma(N)}(z q^{1-N}) \,,
\label{qdet}
\end{eqnarray}
where $\mathfrak{S}_N$ is the set of permutations of $N$ objects.

Conversely, for generic values of the parameters $p$, $q$ and of the central charge $c$, the quantum determinant 
generates the center of the $\Aqp{N}$ algebra.
\end{thm}
The next section is devoted to the proofs of this theorem. 
By generic values of the parameters $p$, $q$ and of the central charge $c$, we mean that there is no functional relation among them, such as the ones used in \cite{AFRS99,AFR} to define deformations of $W_N$ algebras, and that they do not obey any algebraic relation, such as 
 $c=-N$ or $q$ being a root of unity, where it is known that the center is extended, see \cite{FFR} for the former case and \cite{BelaJim,YBS2005} for the latter.

Remark that for $N=2$, it was already proven in \cite{FIJKMY95} that the quantum determinant is central. To the best of our knowledge, the case $N>2$ was not studied yet. 

Taking the limit $p\to 0$, we recover the formula \eqref{qdet-nonell}, as expected.
\begin{coro}\label{coro:deltaqdet}
The quantum determinant is group-like:
\begin{equation}
\Delta \mathrm{qdet} L(z)=\mathrm{qdet} L(zq^{-c^{(2)}/2};pq^{-2c^{(2)}})\otimes \mathrm{qdet} L(zq^{c^{(1)}/2};p),
\end{equation}
where we have specified the $p$-dependence, as in \eqref{coproduct}.
\end{coro}
\prf
We apply the coproduct \eqref{coproduct} to the expression \eqref{qdet0}. For brevity and for this calculation, we note $L(z)\equiv L(z;p)$ and $L^*(z)\equiv L(z;pq^{-2c^{(2)}})$. We get
\begin{eqnarray*}
&&\Delta\mathrm{qdet} L(z)\, A^{(N)}_{N} = \Delta L_1(z)\,\Delta L_2(zq^{-1})\dots \Delta L_N(z q^{1-N})\, A^{(N)}_{N} \\
&& = \Big(1{\otimes}L_1(zq^{c^{(1)}/2})\Big) \Big(L^*_1(zq^{-c^{(2)}/2}){\otimes}1\Big) 
\dots \Big(1{\otimes}L_N(zq^{1-N}q^{c^{(1)}/2})\Big) \Big(L^*_N(zq^{1-N}q^{-c^{(2)}/2}){\otimes}1\Big)\, A^{(N)}_{N} \\
&&\quad =  \Big(1{\otimes}L_1(zq^{-c^{(1)}/2})\dots L_N(z q^{1-N}q^{-c^{(1)}/2})\Big)\Big(L^*_1(zq^{c^{(2)}/2})\dots L^*_N(z q^{1-N}q^{c^{(2)}/2}){\otimes}1\Big)A^{(N)}_{N}
\,  \\
&&\quad =  \Big( 1{\otimes}L_1(zq^{-c^{(1)}/2})L_2(zq^{-1}q^{-c^{(1)}/2})\dots L_N(z q^{1-N}q^{-c^{(1)}/2})A^{(N)}_{N}\Big)\,\Big(\mathrm{qdet} L^*(zq^{c^{(2)}/2}){\otimes}1 \Big)\\
&&\quad =  \Big(1{\otimes}\mathrm{qdet} L(zq^{-c^{(1)}/2})\,A^{(N)}_{N}\Big)\,\Big(\mathrm{qdet} L^*(zq^{c^{(2)}/2}){\otimes}1\Big)
 \\
&&\quad =  \mathrm{qdet} L^*(zq^{-c^{(2)}/2}){\otimes}\,\mathrm{qdet} L(zq^{c^{(1)}/2})\, A^{(N)}_{N}
\end{eqnarray*}
where in  the last step we have used that the quantum determinant is central.
\qed

Theorem \ref{thm:qdet} and corollary \ref{coro:deltaqdet} allow us to introduce the elliptic quantum algebra associated to $\mathfrak{\widehat{sl}}_{N}$:
\begin{defi}
The elliptic quantum algebra $\mathcal{A}_{q,p}(\mathfrak{\widehat{sl}}_{N})$ is the quasi-Hopf algebra defined by the coset
\begin{equation}
\mathcal{A}_{q,p}(\mathfrak{\widehat{sl}}_{N}) = \mathcal{A}_{q,p}(\mathfrak{\widehat{gl}}_{N}) / 
\langle \mathrm{qdet}\,L(z) - q^{c/2} \rangle.
\label{eq:AqpslN}
\end{equation}
\end{defi}
\prf
Since the quantum determinant and $c$ are central, the coset defines an algebra. 
Moreover, $\mathrm{qdet}\,L(z)$ and $q^{c/2}$ are both group-like, so that the coset is a quasi-Hopf algebra.
\qed


\section{Proofs\label{sec:proofs}}
The proof of the main theorem relies on different lemmas and properties.

\subsection{Preliminary lemmas}
\begin{lem}
Let $A_2^{(N)}$ be the antisymmetrizer on $(\CC^N)^2$. Then for generic values of $q$, $\ker \widehat R_{12}(q) = \im A_2^{(N)}$.
\label{lem:ker}
\end{lem}
\prf
It has been shown in \cite{RT86} that $\widehat R_{12}(q)(1-P_{12})=0$, which implies $\ker \widehat R_{12}(q) \supset \im A_2^{(N)}$.
Then, it remains to show that these two spaces have same dimension.
Since the entries of the $R$-matrix of $\mathcal{A}_{q,p}(\mathfrak{\widehat{gl}}_{N})$ are products of (fractional) power of $p$ and analytical functions of $p$, it is sufficient to consider the non-elliptic limit of the matrix, i.e. $p \to 0$, leading to the $R$-matrix \eqref{eq:RUqell} of $\mathcal{U}_{q}(\mathfrak{\widehat{gl}}_{N})$. 

Denoting $\lambda'=\rho_N(q^2)^{-1}\,\lambda$ and $Q=\frac{q}{1+q^2}$, one immediately gets that 
\begin{equation}
\begin{aligned}
\frac{\det(\textsf{R}'(q)-\lambda\id)}{\rho_N(q^2)^{N^2}} &= (1-\lambda')^N \prod_{1 \le i<j \le N}  \big( (Q q^{(2i-2j+N)/N}-\lambda')(Q q^{-(2i-2j+N)/N}-\lambda') - Q^2 \big)
\\
&= (\lambda')^{\frac{N(N-1)}2} (1-\lambda')^N \prod_{1 \le i<j \le N} \big( \lambda'-Q (q^{(2i-2j+N)/N}+q^{-(2i-2j+N)/N}) \big)
\end{aligned}
\end{equation}
from which it follows that for generic values of $q$, the eigenvalue 0 has multiplicity $\half N(N-1)$.
It shows that $\ker \widehat R_{12}(q)$ and $\im A_2^{(N)}$ have same dimension, hence the result.\qed

\begin{lem}
\label{lem:LL}
In the $\Aqp{N}$ algebra, the following identity holds
\begin{equation}
L_{i,j}(z)L_{k,l}(\frac{z}{q})-L_{i,l}(z)L_{k,j}(\frac{z}{q})=L_{k,l}(z)L_{i,j}(\frac{z}{q})-L_{k,j}(z)L_{i,l}(\frac{z}{q}) \quad \forall \,i,j,k,l=1,\dots, N.
\label{proofcondition2}
\end{equation}
In particular, we have 
\begin{equation}
L_{i,j}(z)L_{i,l}(\frac{z}{q})=L_{i,l}(z)L_{i,j}(\frac{z}{q}) \quad \forall \, i,j,l=1,\dots, N.
\label{proofcondition1}
\end{equation}
\end{lem}
\prf 
We consider the RLL relations \eqref{RLL} for $w=z/q$ and project them onto an arbitrary element $e_{i,j}\otimes e_{k,l}$. 
This leads to the following equation, valid for all $i,j,k,l=1,\dots,N$:
\begin{equation}
\sum_{n,m=1}^N \wh{R}_{i,k}^{n,m}(q) L_{n,j}(z) L_{m,l}(\frac{z}{q}) = \sum_{n,m=1}^N \wh{R}^*{}_{m,n}^{j,l}(q)  L_{k,n}(\frac{z}{q})L_{i,m}(z),
\label{RLLPointwise}
\end{equation}
We will refer to this equation as $X_{i,k}^{j,l}(z)$. \\
Note that the lemma \ref{lem:ker} implies the relations 
\begin{equation}
\wh{R}_{i,k}^{j,l}(q) = \wh{R}_{i,k}^{l,j}(q) \,.
\label{MatEntProp}
\end{equation}
Hence, looking at the difference $X_{i,k}^{j,l}(z)-X_{i,k}^{l,j}(z)$, the R.H.S. is equal to
\begin{equation}
\sum_{n,m=1}^N \left({\wh{R}^*}{}_{m,n}^{j,l}(q) - {\wh{R}^*}{}_{m,n}^{l,j}(q)  \right)L_{k,n}(\frac{z}{q})L_{i,m}(z)=0
\end{equation}
In other words, $X_{i,k}^{j,l}(z)-X_{i,k}^{l,j}(z)$ does \textbf{not} depend on $p^*$. Finally, the L.H.S. gives us
\begin{equation}
\sum_{n,m=1}^N \wh{R}_{i,k}^{n,m}(q) \left(L_{n,j}(z) L_{m,l}(\frac{z}{q})-L_{n,l}(z) L_{m,j}(\frac{z}{q})\right)=0.
\label{GaussSystem}
\end{equation}
Note that the indices $j$ and $l$ do not play any role in these relations, 
so if we can solve \eqref{GaussSystem} for one pair $j,l$, we can do it for any. 
We thus consider the equations for  fixed indices $j$ and $l$, and omit them
to ease the notation.

Denoting
\begin{equation}
T_{n,m}(z)\equiv T_{n,m}^{j,l}(z)=L_{n,j}(z) L_{m,l}(\frac{z}{q})-L_{n,l}(z) L_{m,j}(\frac{z}{q}) + (n \leftrightarrow m),
\end{equation}
one gets for fixed $i,j,k,l$, using again the property \eqref{MatEntProp},
\begin{equation}
\sum_{n,m=1}^N \wh{R}_{i,k}^{n,m}(q) T_{n,m}(z)=0.
\label{GaussSystem2}
\end{equation}
Since $T(z)$ is in the symmetric part of $\CC^N\otimes\CC^N$, lemma \ref{lem:ker} implies that the only solution of the linear system \eqref{GaussSystem2} is $T_{n,m}(z)=0$, which is relation \eqref{proofcondition2}.
\qed

\subsection{Explicit expression of the quantum determinant}

The antisymmetrizer $A_N^{(N)}$ in $(\CC^N)^{\otimes N}$ is a rank 1 projector, the eigenvector corresponding to the eigenvalue 1 being given by
\begin{equation}
\sw=\sum_{\sigma \in \mathfrak{S}_N} \text{sgn} (\sigma) e_{\sigma(1)} \otimes \dots \otimes e_{\sigma(N)},
\end{equation}
where $\{e_i\}$ denotes the standard vector basis of $\CC^N$. 
$A_N^{(N)}$ projects any given vector $\sv$ on $\sw$:
\begin{equation}
A_N \sv = \langle \sw,\sv\rangle\, \sw \hspace{1cm} \forall \sv \in (\mathbb{C}^N)^{\otimes N},
\label{AS_Eigenvector}
\end{equation}
where $\langle \sw,\sv\rangle
=\sum_{\sigma \in \mathfrak{S}_N} \text{sgn} (\sigma) v_{\sigma(1) \dots \sigma(N)}$.

Due to  equality \eqref{AS_Eigenvector}, to get an expression for the quantum determinant, it is enough to compute
\begin{align}
L_1(z)\dots L_N(z q^{1-N})\sw=\sum_{i_1,\dots,i_N=1}^N \sum_{\sigma \in \mathfrak{S}_N} \text{sgn} (\sigma) L_{i_1,\sigma(1)}(z) \dots  L_{i_N,\sigma(N)}(q^{1-N}z)(e_{i_1}\otimes \dots \otimes e_{i_N}).
\label{ActionOnEigenvector}
\end{align}
We first prove that all the indices $i_1,\dots,i_N$ in \eqref{ActionOnEigenvector} must be different. 
For such a purpose, we prove that terms with identical indices vanish. The proof is done by recursion on the 'distance' between two identical indices.

Consider the terms with $i_k=i_{k+1}$. Without loss of generality, we can check what happens for $k=N-1$, the reasoning naturally translates to all other possible pairs of adjacent indices. Focusing on the coefficient of $e_{i_1}\otimes \dots \otimes e_{i_{N-2}}\otimes e_{i_N}\otimes e_{i_N}$ only, we write (all indices arbitrary, but fixed)
\begin{align}
\sum_{\sigma \in \mathfrak{S}_N} \text{sgn} & (\sigma) L_{i_1,\sigma(1)}(z) \dots  L_{i_N,\sigma(N-1)}(q^{2-N}z) L_{i_N,\sigma(N)}(q^{1-N}z)=\nonumber\\
=\ &\sum_{\sigma' \in \mathfrak{S}_N} \text{sgn} (\sigma' \circ s_{N,N-1}) L_{i_1,\sigma'(1)}(z) \dots  L_{i_N,\sigma'(N)}(q^{2-N}z) L_{i_N,\sigma'(N-1)}(q^{1-N}z)\nonu
=\ &\frac{1}{2}\sum_{\sigma' \in \mathfrak{S}_N} \text{sgn} (\sigma') L_{i_1,\sigma'(1)}(z) \dots  L_{i_{N-2},\sigma'(N-2)}(q^{3-N}z)\times\nonumber\\
&\times \Big( L_{i_N,\sigma'(N-1)}(q^{2-N}z)  L_{i_N,\sigma'(N)}(q^{1-N}z)- L_{i_N,\sigma'(N)}(q^{2-N}z)  L_{i_N,\sigma'(N-1)}(q^{1-N}z) \Big)\nonu
=\ &0
\label{Case1AntiSym}
\end{align}
where the last equality is done by virtue of \eqref{proofcondition1}.

Suppose now that the terms where $i_k=i_{k+n}$ with $1\leq n\leq m$ have zero contribution and consider the term where $i_k=i_{k+m+1}$:
\begin{align}
\sum_{\sigma \in \mathfrak{S}_N} \text{sgn} & (\sigma) L_{i_1,\sigma(1)}(z) \dots  L_{i_{k},\sigma(k)}(q^{1-k}z)\dots L_{i_{k},\sigma(k+m+1)}(q^{-m-k}z)\dots L_{i_N,\sigma(N)}(q^{1-N}z) \nonumber \\
=\ &\frac{1}{2}\sum_{\sigma' \in \mathfrak{S}_N} \text{sgn} (\sigma') L_{i_1,\sigma'(1)}(z) \dots  L_{i_{k-1},\sigma'(k-1)}(q^{2-k}z) \nonumber \\
&\times \Big( L_{i_k,\sigma'(k)}(q^{1-k}z)  L_{i_{k+1},\sigma'(k+1)}(q^{-k}z)- L_{i_k,\sigma'(k+1)}(q^{1-k}z)  L_{i_{k+1},\sigma'(k)}(q^{-k}z) \Big) \nonu
&\times  L_{i_{k+2},\sigma'(k+2)}(q^{3-k}z)\dots L_{i_{k},\sigma'(k+m+1)}(q^{-m-k}z)\dots L_{i_N,\sigma'(N)}(q^{1-N}z)\label{toto} \\
=\ &-\frac{1}{2}\sum_{\sigma' \in \mathfrak{S}_N} \text{sgn} (\sigma') L_{i_1,\sigma'(1)}(z) \dots  L_{i_{k-1},\sigma'(k-1)}(q^{2-k}z) \nonumber \\
&\times \Big( L_{i_{k+1},\sigma'(k)}(q^{1-k}z)  L_{i_{k},\sigma'(k+1)}(q^{-k}z)- L_{i_{k+1},\sigma'(k+1)}(q^{1-k}z)  L_{i_{k},\sigma'(k)}(q^{-k}z) \Big)\nonu
&\times  L_{i_{k+2},\sigma'(k+2)}(q^{3-k}z)\dots L_{i_{k},\sigma'(k+m+1)}(q^{-m-k}z)\dots L_{i_N,\sigma'(N)}(q^{1-N}z)\label{titi} \\
=\ &-\sum_{\sigma \in \mathfrak{S}_N} \text{sgn} (\sigma) L_{i'_1,\sigma(1)}(z) \dots  L_{i'_{k},\sigma(k)}(q^{1-k}z)\dots L_{i'_{k},\sigma(k+m)}(q^{-m-k}z)\dots L_{i'_N,\sigma(N)}(q^{1-N}z) \,.
\end{align}
To get \eqref{toto}, we have used the same trick as in the calculation of \eqref{Case1AntiSym}.
To go from \eqref{toto} to \eqref{titi} we have used the relation \eqref{proofcondition2}. 
In the last equality, we have introduced the indices
$i'_\ell=i_\ell$ for $\ell\notin \{k,k+1\}$ and $i'_k=i_{k+1}$, $i'_{k+1}=i_{k}$. This last expression vanishes due to the recursion hypothesis.

Since all indices $i_r$ are different, we can replace the sum on $i_1,\dots,i_N$ by a sum over permutations 
 $\mu \in \mathfrak{S}_N$. We pick one such permutation and examine the coefficient of 
 $e_{\mu(1)}\otimes \dots \otimes e_{\mu(N)}$:
\begin{align}
\chi_{\mu} :=\ & \sum_{\sigma \in \mathfrak{S}_N} \text{sgn} (\sigma)L_{\mu(1),\sigma(1)} \dots L_{\mu(N),\sigma(N)}\\
=\ &\frac{1}{2}\sum_{\sigma \in \mathfrak{S}_N} \text{sgn} (\sigma)L_{\mu(1),\sigma(1)} \dots L_{\mu(k-1),\sigma(k-1)} 
\nonumber \\
& \quad \times \big\{ L_{\mu(k),\sigma(k)}L_{\mu(k+1),\sigma(k+1)}-L_{\mu(k),\sigma(k+1)}L_{\mu(k+1),\sigma(k)}\big\}
L_{\mu(k+2),\sigma(k+2)}\dots L_{\mu(N),\sigma(N)} \,.
\end{align}
But we can also look at a different permutation $\mu'=\mu \circ s_k$. In this case, we find that
\begin{align}
\chi_{(\mu \circ s_k)} =\ &\frac{1}{2}\sum_{\sigma \in \mathfrak{S}_N} \text{sgn} (\sigma)L_{\mu(1),\sigma(1)} \dots L_{\mu(k-1),\sigma(k-1)} \nonumber \\
& \quad \times \big\{ L_{\mu(k+1),\sigma(k)}L_{\mu(k),\sigma(k+1)}-L_{\mu(k+1),\sigma(k+1)} L_{\mu(k),\sigma(k)}\big\}
L_{\mu(k+2),\sigma(k+2)}\dots L_{\mu(N),\sigma(N)} \,.
\end{align}
Once more, condition \eqref{proofcondition2} shows that $\chi_{(\mu \circ s_k)}=-\chi_{\mu}$. This allows us to conclude that in fact, for any $\sigma, \mu \in \mathfrak{S}_N$, we have $\chi_{\mu}=\text{sgn}(\sigma) \chi_{(\mu \circ \sigma)}$. In particular, $\chi_{\mu}=\text{sgn}(\mu) \chi_{\text{id}}$, and we finally arrive at the following result:
\begin{align}
L_1(z)\dots L_N(z q^{1-N})\sw &= \frac{1}{N!}\,\sum_{\mu \in \mathfrak{S}_N}  \chi_{\mu}\,
e_{\mu(1)}\otimes \dots \otimes e_{\mu(N)} \\
&= \frac{1}{N!}\,\sum_{\mu \in \mathfrak{S}_N}  \text{sgn}(\mu)\, \chi_{\text{id}}\,
e_{\mu(1)}\otimes \dots \otimes e_{\mu(N)}= \chi_{\text{id}}\, \sw \,.
\end{align}
From this, we directly infer that the quantum determinant is $\chi_{\text{id}}$, which proves the equality \eqref{qdet}. \qed

Remark that in this way we have proved that  
\begin{equation}\label{qdet-scalar}
L_1(z)\dots L_N(z q^{1-N})\, A^{(N)}_{1\dots N} = \mathbb{M}(z)\, A^{(N)}_{1\dots N}\,,
\end{equation}
where $\mathbb{M}(z)$ is scalar in the spaces 1,...,$N$ and given by  \eqref{qdet}. 
It remains to prove that $\mathbb{M}(z)$ is central in $\Aqp{N}$: it is done in section \ref{sect:central}.

\subsection{Value of q-det$L(u)$ in the fundamental representation}
\begin{lem}\label{lem:RAntisym}
The $R$-matrix \eqref{RMatrixAN} for the $\Aqp{N}$ algebra obeys the following relation
\begin{equation}
\wh{R}_{10}\left(z\right) \dots\wh{R}_{N0}\left(zq^{1-N}\right) A^{(N)}_{1\dots N}
=A^{(N)}_{1\dots N} \,.
\label{AntisymProp}
\end{equation}
\end{lem}
\prf 
We apply the evaluation maps $\pi_j$: $L_j(z)\to\wh{R}_{j0}(z),$ $j=1,...,N$ to the equality \eqref{qdet-scalar}:
\begin{equation}\label{eval-qdet}
\wh{R}_{10}\left(z\right) \dots\wh{R}_{N0}\left(zq^{1-N}\right) A^{(N)}_{1\dots N}
=\pi\big(\text{qdet}L(z)\big)\,A^{(N)}_{1\dots N}=M_0(z)\,A^{(N)}_{1\dots N} \,,
\end{equation}
where $\pi=\pi_1\otimes\pi_2\otimes...\otimes\pi_N$ and $M_0(z)$ is a matrix $M(z)$ (yet to be determined)
acting on the space 0 only. One can show from the evaluation of the relation \eqref{qdet} that $M(z)$ is a diagonal matrix:
\begin{align}
M(z) &= \prod_{j=0}^{N-1}\eta\left(\frac{z}{q^j}\right) \sum_{k=1}^N q^{2k-N-1} \sum_{\sigma \in \mathfrak{S}N} 
\text{sgn} (\sigma)\wh{S}_{1,k}^{\sigma(1)}(z)\dots \wh{S}_{N,k+\sum_{i=1}^{N-1} (i-\sigma(i))}^{\sigma(N)} 
\left(\frac{z}{q^{N-1}}\right) e_{k,k} \nonumber \\
&\equiv \sum_{k=1}^N m_k(z)\,e_{k,k} \,,
\end{align}
where we set, for \emph{any} indices $a,b,c$,
\begin{align}
\wh{S}_{a,c}^b(z):=\frac{\Theta_{p^N}(p^{N+c-a} q^2 z^2)}{\Theta_{p^N}(p^{N+c-b} z^2)\Theta_{p^N}(p^{N+b-a} q^2)} \,.
\end{align}
Using the expression of $\eta(z)$, we get
\begin{equation}\label{expression:m}
m_k(z)=\left(\!-\frac{(p^N;p^N)_\infty}{(p;p)_\infty}\right)^{3N}\,q^{2k-2N}\,\frac{\Theta_p(q^2)^N\,\Theta_p(z^2)}{\Theta_p(q^2z^2)}
\sum_{\sigma \in \mathfrak{S}N} \text{sgn}(\sigma) \prod_{\ell=1}^N \wh{S}_{\ell,k+\sum_{i=1}^{\ell-1} (i-\sigma(i))}^{\sigma(\ell)}(\frac{z}{q^{\ell-1}}) \,.
\end{equation}
Now, from the invariance property \eqref{eq:h-inv} of the $R$-matrix, it is easy to show that
\begin{align}
h_0h_1\cdots h_N \,\pi\big(\text{qdet}L(z)\big)\,A^{(N)}_{1\dots N}
& = \pi\big(\text{qdet}L(z)\big)\,h_0h_1\cdots h_N \,A^{(N)}_{1\dots N} \nonumber \\
& = \pi\big(\text{qdet}L(z)\big) \,A^{(N)}_{1\dots N}\,h_0h_1\cdots h_N \,.
\end{align}
Due to the expression \eqref{eval-qdet}, it implies $[h_0\,,\,M_0(z)]=0$ which can be recasted as $m_k(z)=m_{k+1}(z)$, that is to say $M(z)=m(z)\,\II$. 

Using property \eqref{id:theta}, it is easy to show that 
$m(z;qp^{N/2},p)=m(z;q,p)$, where we have indicated explicitly the dependence in the parameters $q$ and $p$.
In other words, since $|p|<1$, we have
\begin{equation}\label{limq0}
m(z;q,p)=m(z;qp^{\ell N/2},p)=\lim_{\ell\to\infty}m(z;qp^{\ell N/2},p)=\lim_{q'\to0}m(z;q',p),\quad\forall\ p,q.
\end{equation}
This shows that $m(z;q,p)$ does not depend on $q$. To compute it, we take the limit $q\to1$. 
To show that this limit is well-defined, we computed explicitly 
relation \eqref{AntisymProp} in the limit $p\to 0$ and generic values of $q$ and $N$ (it corresponds to the 
non-elliptic presentation of $\Uq{N}$, see section \ref{sect:non-ell}). We got $m(z;q,p)|_{p=0}=1$. Since $p$, $q$ and $c$ are generic, it
shows that the limit \eqref{limq0} exists at least in a neighborhood of $p=0$.
Due to the term 
$\Theta_p(q^2)^N$ in \eqref{expression:m}, which vanishes in the limit $q\to1$, one sees that only the term $\sigma=id$ contributes to $m(z;q,p)$.
Then, a direct calculation shows that $m(z;1,p)=1$ for generic values of $p$ and $N$.
\qed

We checked relation \eqref{AntisymProp} for $N=2,3$ and generic values of $p$ and $q$.

\subsection{Centrality of the quantum determinant\label{sect:central}}
We wish to show that, for $z,w \in \mathbb{C}$ 
\begin{equation}
\left[\mathrm{qdet}L(z), L_{0}(w)\right]=0 \,.
\label{commutator}
\end{equation}
This will be achieved by commuting $L_0(w)$ through the expression \eqref{qdet2} for the quantum determinant:
\begin{align}
\mathrm{qdet}L(z)\, & L_{0}(w) \;=\; 
\mathrm{tr}_{1\dots N}\left[L_1(z)\dots L_N(z q^{1-N}) L_0(w) A^{(N)}_{1\dots N}\right] 
\nonumber \\
=\ &\mathrm{tr}_{1\dots N}\left[L_1(z)\dots L_{N-1}(z q^{2-N})\wh{R}^{-1}_{N0}\left(\frac{z}{w }q^{1-N}\right) L_0(w) L_N(z q^{1-N})\wh{R}^{*}_{N0}\left(\frac{z}{w }q^{1-N}\right) A^{(N)}_{1\dots N}\right]
\nonumber \\
=\ &\mathrm{tr}_{1\dots N}\left[\wh{R}^{-1}_{N0}\left(\frac{z}{w }q^{1-N}\right)\dots\wh{R}^{-1}_{10}\left(\frac{z}{w}\right)L_0(w) L_1(z)\dots L_N(z q^{1-N})\right.
\nonumber \\
&\quad\quad\quad \times\left.
\wh{R}^{*}_{10}\left(\frac{z}{w }\right) \dots\wh{R}^{*}_{N0}\left(\frac{z}{w }q^{1-N}\right) A^{(N)}_{1\dots N}\right]
\nonumber \\
=\ &\mathrm{tr}_{1\dots N}\left[\wh{R}^{-1}_{N0}\left(\frac{z}{w }q^{1-N}\right)\dots\wh{R}^{-1}_{10}\left(\frac{z}{w}\right)L_0(w) L_1(z)\dots L_N(z q^{1-N}) \, A^{(N)}_{1\dots N}\right],
\nonumber \\
=\ &\mathrm{tr}_{1\dots N}\left[\wh{R}^{-1}_{N0}\left(\frac{z}{w }q^{1-N}\right)\dots\wh{R}^{-1}_{10}\left(\frac{z}{w}\right)L_0(w) \,\mathrm{qdet}L(z) \, A^{(N)}_{1\dots N}\right],
\end{align}
where we used the RLL relations \eqref{RLL} and the fact that generators acting in different subspaces commute. The last equalities are due to lemma \ref{lem:RAntisym} and definition \ref{qdet0}. \\
Next, using the fact that the quantum determinant is a scalar in the spaces $0,1,\dots,N$, we get
\begin{align}
\mathrm{qdet}L(z)\, L_{0}(w)=\ &\mathrm{tr}_{1\dots N}\left[ \wh{R}^{-1}_{N0}\left(\frac{z}{w }q^{1-N}\right)\dots
\wh{R}^{-1}_{10}\left(\frac{z}{w}\right)A^{(N)}_{1\dots N} L_0(w)\right] \mathrm{qdet}L(z)
\nonumber \\
=\ &L_0(w)\mathrm{tr}_{1\dots N}\left[A^{(N)}_{1\dots N}\right] \mathrm{qdet}L(z)
\nonumber \\
=\ &L_0(w)\, \mathrm{qdet}L(z),
\label{final}
\end{align}
where we used that $L_0(w)$ and the antisymmetrizer commute as they live in different spaces, applied the inverse of \eqref{AntisymProp}, and finally traced over the antisymmetrizer. As a consequence, 
 the quantum determinant lies in the center of the algebra $\Aqp{N}$ as desired.

\subsection{Center of the algebra $\Aqp{N}$}
It is known that in $\Uq{N}$ and for generic values of $q$ and $c$, the quantum determinant, as defined in \eqref{qdet-homo}, generates the center of this quantum algebra \cite{FRT}. Since, being a twist of it, $\Aqp{N}$ is isomorphic to $\Uq{N}$ as an algebra \cite{JKOS} and for generic values of $p$, $q$ and $c$ (in the sense explained in section \ref{sec:main}), \eqref{qdet-homo} also describes the full center of the algebra $\Aqp{N}$. The same is true for the other two expressions \eqref{qdet-princ} and \eqref{qdet-nonell}, that are just the same quantum determinant in different presentations.

Moreover, we have shown that expression \eqref{qdet} also lies in the center of $\Aqp{N}$, and that its limit for $p\to 0$ coincides with the quantum determinant of $\Uq{N}$ in the non-elliptic presentation, expression \eqref{qdet-nonell}. It means that the limit $p\to 0$ defines a surjective mapping from a set of elements (defined by the expression \eqref{qdet}) in the center of $\Aqp{N}$ to a generating set (defined by \eqref{qdet-nonell}) of the same center. 

Thus, it implies that \eqref{qdet} also defines a generating set of the center of $\Aqp{N}$ for generic values of $p$, $q$ and $c$.


\end{document}